# Système dérivé et suite duale d'une suite barypolygonale - Partie 1

Par David Pouvreau[1] et Vincent Bouis[2]

**Résumé**

La présente étude prolonge trois récents articles ayant défini les suites barypolygonales et établi leurs propriétés de convergence. Une suite barypolygonale $\mathfrak{B}$ quelconque d'un ensemble fini $\mathcal{A}$ de $p \geq 2$ points d'un espace affine de dimension finie quelconque étant donnée, on peut définir par récurrence une certaine suite de suites barypolygonales initialisée en $\mathfrak{B}$. Cette suite $\left(\mathfrak{B}^{(m)}\right)_{m \in \mathbb{N}}$, appelée suite des dérivées de $\mathfrak{B}$, est déterminée par des suites de réels solutions d'un système récurrent non linéaire $(S)$ : le système barypolygonal dérivé de $\mathfrak{B}$. Chaque terme de la suite $\left(\mathfrak{B}^{(m)}\right)_{m \in \mathbb{N}}$ converge vers un point $G_m$. La suite $(G_m)_{m \in \mathbb{N}}$ est appelée la suite duale de $\mathfrak{B}$. La convergence de cette dernière suite et les propriétés du système dérivé sont ici étudiées pour tout $p$ si $\mathfrak{B}$ est régulière et dans tous les cas où $p \in \{2; 3\}$.

**Mots clefs** : Suites barypolygonales, Systèmes dynamiques discrets

**Abstract**

This study continues three recent papers in which barypolygonal sequences have been defined and their properties of convergence demonstrated. Any barypolygonal sequence $\mathcal{B}$ of a finite set $\mathcal{A}$ comprising $p \geq 2$ points of any finite dimensional affine space can be used in order to define recurrently a definite sequence of barypolygonal sequences starting with $\mathfrak{B}$. This sequence $\left(\mathfrak{B}^{(m)}\right)_{m \in \mathbb{N}}$, called sequence of $\mathfrak{B}$'s derivatives, is determined by real sequences that are solutions of a non linear recurrent system $(S)$ : the barypolygonal derived system of $\mathfrak{B}$. Each term of the sequence $\left(\mathfrak{B}^{(m)}\right)_{m \in \mathbb{N}}$ converges toward a point $G_m$. The sequence $(G_m)_{m \in \mathbb{N}}$ is the dual sequence of $\mathfrak{B}$. The convergence of the latter and the properties of the derived system are here investigated for any $p$ if $\mathfrak{B}$ is regular and in any case if $p \in \{2; 3\}$.

## 1. Introduction

Cet article prolonge les trois études récemment publiées qui concernent les suites barypolygonales : (Pouvreau, 2016), (Pouvreau, Eupherte, 2016) et (Bouis, 2017). Il convient au préalable de renouveler la présentation des objets qui seront considérés ici.

Soit $E$ un espace affine réel de dimension finie quelconque. Soit $p \in \mathbb{N} \backslash \{0\}$. On considère une famille ordonnée $\mathcal{A}$ de points distincts $(A_k)_{1 \leq k \leq p}$ de $E$ et une famille ordonnée $t = (t_k)_{1 \leq k \leq p}$ de réels de $]0; 1[$. On appelle « $t$–barypolygone » de $\mathcal{A}$ la famille ordonnée de points $\mathcal{B} = (B_k)_{1 \leq k \leq p}$ dont les éléments sont définis par les conditions barycentriques :

$$\begin{cases} \forall\, k \in [\![1; p-1]\!], \ B_k = \mathrm{bar}\,\{(A_k\,;t_k)\,;(A_{k+1}\,;1-t_k)\} \\ \qquad\qquad B_p = \mathrm{bar}\,\{(A_p\,;t_p)\,;(A_1\,;1-t_p)\} \end{cases}$$

On appelle alors suite « $t$–barypolygonale » de $\mathcal{A}$ la suite $\mathfrak{B} = (\mathcal{B}_n)_{n \in \mathbb{N}}$ de familles ordonnées de points définie par récurrence de la manière suivante :

---

[1] Professeur agrégé de mathématiques et docteur en histoire des sciences. Centre Universitaire de Mayotte et Institut Alexander Grothendieck (Université de Montpellier). Email : david_pouvreau@orange.fr
[2] Etudiant à l'Ecole Normale Supérieure de Paris. Email : vbouis@clipper.ens.fr

$$\begin{cases} \mathcal{B}_0 = \mathcal{A} \\ \forall\, n \in \mathbb{N}, \ \mathcal{B}_{n+1} \text{ est le } t\text{– barypolygone de } \mathcal{B}_n \end{cases}$$

Cette suite $\mathfrak{B}$ est dite régulière si $t_i = t_j$ pour tout $(i;j) \in [\![1;p]\!]^2$. On convient alors de noter $t$ la valeur commune des paramètres $(t_k)_{1 \leq k \leq p}$. On dit que $\mathfrak{B}$ est irrégulière sinon.

Le théorème suivant (dit « de convergence barypolygonale ») a été démontré de deux manières dans (Pouvreau, Eupherte, 2016) et d'une troisième manière dans (Bouis, 2017)[3] :

**Théorème 1**

Pour toute famille $\mathcal{A}$ de points distincts et toute famille ordonnée $t = (t_k)_{1 \leq k \leq p}$ de réels de $]0;1[$, la suite « $t$–barypolygonale » $\mathfrak{B}$ de $\mathcal{A}$ converge vers le point

$$G = \mathrm{bar}\left\{\left(A_k; \frac{1}{1-t_k}\right)\right\}_{1 \leq k \leq p} = \mathrm{bar}\left\{\left(A_k; \prod_{\substack{1 \leq i \leq p \\ i \neq k}} (1-t_i)\right)\right\}_{1 \leq k \leq p}$$

La figure 1 illustre ce théorème en dimension 2 dans le cas où $p = 5$, avec une suite $\left(\frac{1}{61}; \frac{1}{41}; \frac{1}{28}; \frac{1}{19}; \frac{1}{13}\right)$-barypolygonale de $\mathcal{A}$.

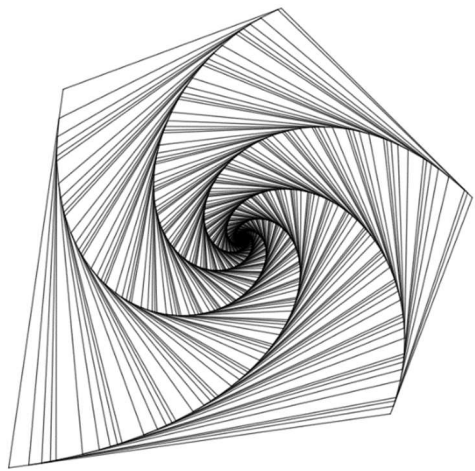

Figure 1

Cet article vise à mettre en évidence une propriété remarquable de « dualité » engendrée par toute suite $t$-barypolygonale. Commençons par préciser ce qui est entendu de la sorte.

## 2. Système dérivé, et duale d'une suite barypolygonale : définitions

On suppose que $p \geq 2$. Notons alors :

$$\forall\, k \in [\![1;p]\!], \ t_k^{(1)} = \prod_{\substack{1 \leq i \leq p \\ i \neq k}} (1-t_i)$$

On a : $\forall\, k \in [\![1;p]\!], \ t_k^{(1)} \in ]0;1[$. Notons $t^{(1)}$ la famille ordonnée $t^{(1)} = \left(t_k^{(1)}\right)_{1 \leq k \leq p}$. On peut alors aussi définir la suite « $t^{(1)}$-barypolygonale » de $\mathcal{A}$, dite « *suite dérivée première* » de $\mathfrak{B}$ et notée $\mathfrak{B}^{(1)}$.

---

[3] Ce théorème ayant de plus été généralisé dans (Bouis, 2017) au cas des suites $k$-barypolygonales avec $k \geq 2$, pour lesquelles le nombre de points consécutifs considérés dans les systèmes barycentriques est $k$ : le présent article est exclusivement consacré aux suites 2-barypolygonales en ce sens.

Soit $\left(t^{(m)}\right)_{m\in\mathbb{N}}$ la suite de $(]0;1[^p)^{\mathbb{N}}$ définie par :

$$(S) : \begin{cases} t^{(0)} = t \\ \forall\, m \in \mathbb{N},\ t^{(m+1)} = \left(t_k^{(m+1)}\right)_{1\leq k\leq p} \quad \text{où } \forall\, k \in [\![1;p]\!],\ t_k^{(m+1)} = \prod_{\substack{1\leq i\leq p \\ i\neq k}} \left(1 - t_i^{(m)}\right) \end{cases}$$

On appellera $(S)$ le « *système barypolygonal dérivé* » de $\mathfrak{B}$ (ou de la famille $t$). Pour tout $m \in \mathbb{N}$, on peut définir au moyen de $t^{(m)}$ la « *suite dérivée m-ième* » de $\mathfrak{B}$, que l'on note $\mathfrak{B}^{(m)}$, par :

Pour tout $m \in \mathbb{N}$, $\mathfrak{B}^{(m)}$ est la suite $t^{(m)}$–barypolygonale de $\mathcal{A}$.

Pour chaque $m \in \mathbb{N}$, la suite $\mathfrak{B}^{(m)}$ a un point limite que l'on notera $G_m$, déterminé par le théorème de convergence barypolygonale. On appelle *suite duale* de $\mathfrak{B}$ la suite $(G_m)_{m\in\mathbb{N}}$.

## 3. Les deux problèmes : énoncés, résultats préliminaires et conjectures

Deux problèmes sont naturellement posés dans le prolongement de ces définitions :

(1) Quelle est la dynamique du système dérivé d'une suite barypolygonale ?
(2) Quel est le comportement asymptotique de sa suite duale ?

Deux cas peuvent d'emblée être résolus. D'abord celui où $p = 2$, dont la preuve est évidente.

**Théorème 2**
Soit $\mathfrak{B}$ une suite $t$-barypolygonale de $\mathcal{A}$ avec $\text{Card}(\mathcal{A}) = 2$ et $t = (t_1; t_2)$.
Si $t_2 = 1 - t_1$, le système dérivé et la duale de la suite $\mathfrak{B}$ sont stationnaires.
Si $t_2 \neq 1 - t_1$, le système dérivé et la duale de la suite $\mathfrak{B}$ sont 2-périodiques.

Le cas des suites $t$-barypolygonale *régulières* lorsque $p \geq 3$ peut également être d'ores et déjà considéré. La démonstration du théorème suivant s'avérera par la suite utile dans l'étude du cas général :

**Théorème 3**
Soit $\mathfrak{B}$ une suite $t$-barypolygonale régulière de $\mathcal{A}$, avec $\text{Card}(\mathcal{A}) = p \geq 3$. Alors la suite duale de $\mathfrak{B}$ converge vers le centre de gravité $G$ de $\mathcal{A}$. De plus, $\alpha_p$ constituant l'unique solution dans $[0;1]$ de l'équation $(E_p) : x^{p-1} + x - 1 = 0$, le système barypolygonal dérivé de $\mathfrak{B}$ :
(i)    est stationnaire si $t = 1 - \alpha_p$, au vecteur de $\mathbb{R}^p$ de coordonnées égales à $1 - \alpha_p$ ;
(ii)   diverge si $t \neq 1 - \alpha_p$, les $p$ suites extraites identiques $\left(t_k^{(2m)}\right)_{1\leq k\leq p}$ convergeant avec la même monotonie vers 0 ou 1, et les $p$ suites extraites identiques $\left(t_k^{(2m+1)}\right)_{1\leq k\leq p}$ convergeant vers 1 ou 0 et avec une variation contraire de celle des suites $\left(t_k^{(2m)}\right)_{1\leq k\leq p}$.

<u>Démonstration.</u>
Compte tenu de l'identité de tous les $t_k^{(m)}$ pour chaque $m \in \mathbb{N}$ et du théorème de convergence barypolygonale, il est immédiat que la suite duale est stationnaire, son lieu étant le point $G$.

Considérons maintenant le système dérivé. Notons, en tenant compte de l'identité des $t_k$ :

$$\forall\, k \in [\![1;p]\!],\ \forall\, m \in \mathbb{N},\ u_m = 1 - t_k^{(m)}$$

Alors $u_0 = 1 - t$ et : $\forall\, m \in \mathbb{N},\ u_{m+1} = 1 - u_m^{p-1}$. Soit $f_p$ définie sur $[0;1]$ par : $f_p(x) = 1 - x^{p-1}$. C'est une bijection décroissante de $[0;1]$ sur lui-même. On a : $\forall\, m \in \mathbb{N},\ u_{m+1} = f_p(u_m)$. On montre facilement que la fonction $x \mapsto f_p(x) - x$ est une bijection de $[0;1]$ sur $[-1;1]$. Donc $f_p$ admet pour unique point fixe dans $[0;1]$ la solution $\alpha_p$ dans cet intervalle de l'équation $(E_p)$ : $x^{p-1} + x - 1 = 0$. $f_p$ étant continue sur $[0;1]$, la suite $(u_m)$, si elle converge, ne peut donc avoir pour limite que $\alpha_p$.

Si $u_0 = \alpha_p$, la suite $(u_m)$ est stationnaire. Supposons que $u_0 \neq \alpha_p$. La fonction $f_p$ décroît strictement, donc $g_p = f_p \circ f_p$ est une bijection strictement croissante de $[0;1]$ sur $[0;1]$. Il en résulte par récurrence que les suites $(u_{2m})$ et $(u_{2m+1})$ sont monotones. Étant bornées, elles convergent vers une solution dans $[0;1]$ de l'équation $h_p(x) = 0$, où $h_p$ est la fonction définie sur $[0;1]$ par :

$$h_p(x) = g_p(x) - x = 1 - x - (1 - x^{p-1})^{p-1}$$

0, 1 et $\alpha_p$ sont solutions de cette équation. Ce sont les seules. En effet, on montre facilement que $h_p'$ est négative en 0 et en 1, et qu'elle a un maximum en $\theta_p = (1/p)^{1/(p-1)}$ qui vaut :

$$\mu(p) = \left(\frac{p-1}{\frac{p-2}{p^{p-1}}}\right)^p - 1$$

On peut alors établir par étude de variations classique que : $\forall p \in \mathbb{N} \setminus \{0;1;2\}, \ \mu(p) > 0$.

D'où l'existence d'exactement deux annulations de $h_p'$ sur $[0;1]$, l'une en un réel $\beta_p \in ]0;\alpha_p[$, l'autre en un réel $\gamma_p \in ]\alpha_p;1[$ ; de sorte que le tableau de variation de $h_p$ est :

| $x$ | 0 | | $\beta_p$ | | $\alpha_p$ | | $\gamma_p$ | | 1 |
|---|---|---|---|---|---|---|---|---|---|
| $h_p'(x)$ | | $-$ | 0 | $+$ | | $+$ | 0 | $-$ | |
| $h_p$ | 0 ↘ $h_p(\beta_p)$ | | | ↗ 0 ↗ | | | $h_p(\beta_p)$ ↘ | | 0 |

On en déduit le résultat annoncé concernant les annulations de $h_p$, ainsi que son signe.

Si $0 < u_0 < \alpha_p$, alors on montre par récurrence au moyen de la fonction $f_p$ que :

$$\forall m \in \mathbb{N}, \ (0 < u_{2m} < \alpha_p) \ \text{ et } \ (\alpha_p < u_{2m+1} < 1)$$

Il résulte alors de l'étude du signe de $h_p$ que $(u_{2m})$ est strictement décroissante et que $(u_{2m+1})$ est strictement croissante. La seule possibilité étant que la première suite converge vers 0 et que la seconde converge vers 1. Les résultats sont les mêmes à permutation près des rôles si $\alpha_p < u_0 < 1$.

Les figures 2 à 7 illustrent successivement, en dimension 2 et dans un cas où $p = 4$, les dérivées d'ordre 0 à 5 d'une suite barypolygonale régulière d'un ensemble $\mathcal{A}$, le paramètre choisi étant ici 0,2 :

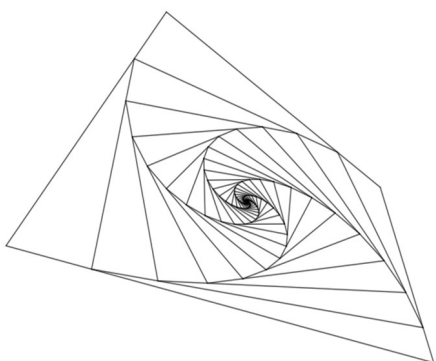 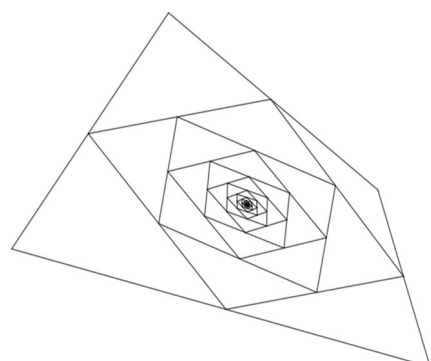

Figure 2                                      Figure 3

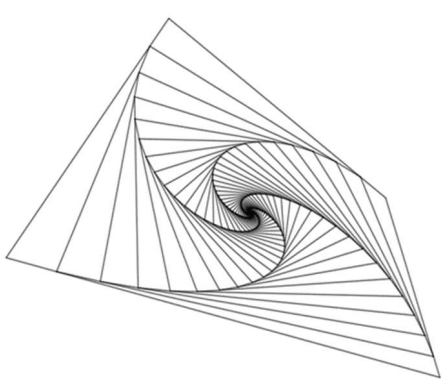
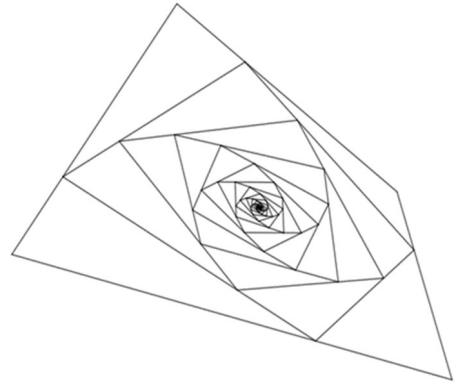

Figure 4                                                            Figure 5

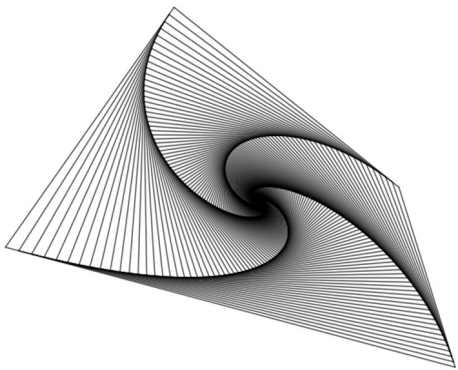
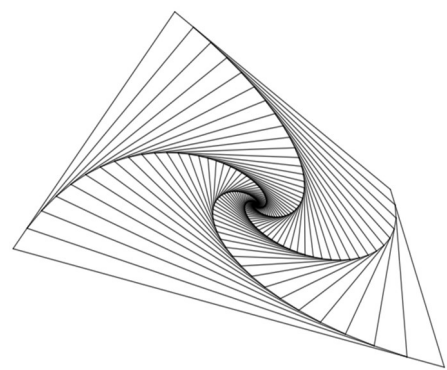

Figure 6                                                            Figure 7

Il découle des considérations précédentes que les problèmes ne demeurent désormais posés que lorsque $p \geq 3$ et dans le cas des suites barypolygonales *irrégulières*. On se placera donc toujours dans cette hypothèse dans le reste de cette étude, sauf mention contraire.

L'examen de quelques cas particuliers s'avère utile pour comprendre dans ces cas le comportement des $p$ suites de coefficients $\left(t_1^{(m)}\right)_{m \in \mathbb{N}}, \ldots, \left(t_p^{(m)}\right)_{m \in \mathbb{N}}$, ainsi que celui de la suite duale. Le lecteur pourra vérifier que ces examens conduisent à la double conjecture suivante, que cet article va mettre à l'épreuve. Si $p \geq 3$, alors :

(1) Les suites $\left(t_k^{(m)}\right)_{m \in \mathbb{N}}$ (où $k \in [\![1;p]\!]$) solutions du système dérivé d'une suite barypolygonale donnée ne convergent que si cette suite est régulière et de paramètre $(1 - \alpha_p)$, où $\alpha_p$ est l'unique solution dans $[0; 1]$ de l'équation $(E_p) : x^{p-1} + x - 1 = 0$. Ces $p$ suites sont alors stationnaires. Sinon, elles ont exactement deux valeurs d'adhérence, 0 et 1, qui sont alternativement les limites des suites extraites des termes de rangs pairs ou des suites extraites des termes de rangs impairs dans les $p$ cas.

(2) La suite duale d'une suite barypolygonale quelconque de $\mathcal{A}$ converge toujours vers le centre de gravité de $\mathcal{A}$.

## 4. Dynamique des systèmes barypolygonaux dérivés et convergence duale lorsque $p = 3$

L'étude détaillée du cas général où $p \geq 3$, assez longue, fera l'objet d'une seconde partie publiée dans un prochain numéro de *Quadrature*. Nous nous limiterons ici à étudier le cas où $p = 3$, introduction indispensable aux calculs assez techniques nécessaires à cette généralisation.

On dispose alors de trois suites $\left(t_1^{(m)}\right)$, $\left(t_2^{(m)}\right)$ et $\left(t_3^{(m)}\right)$ d'éléments de $]0;1[$ qui sont les solutions du système dérivé d'une suite $t$-barypolygonale $\mathfrak{B}$ d'une famille ordonnée de points $\mathcal{A} = (A_k)_{1 \leq k \leq 3}$, où $t$ est un triplet $(t_1; t_2; t_3)$ de $]0;1[^3$ dont les éléments ne sont pas tous égaux.

Notons, pour tout $m \in \mathbb{N}$ :

$$u_m = 1 - t_1^{(m)} \quad ; \quad v_m = 1 - t_2^{(m)} \quad ; \quad w_m = 1 - t_3^{(m)}$$

Les trois suites ainsi déterminées sont les solutions du système récurrent (dit « conjugué » de $(S)$) :

$$(\Sigma) : \begin{cases} u_0 = 1 - t_1 \;\; ; \;\; v_0 = 1 - t_2 \;\; ; \;\; w_0 = 1 - t_3 \\ u_{m+1} = 1 - v_m w_m \\ v_{m+1} = 1 - u_m w_m \\ w_{m+1} = 1 - u_m v_m \end{cases}$$

C'est ce système que nous allons étudier, les conclusions concernant $(S)$ et la suite duale de $\mathfrak{B}$ s'en déduisant immédiatement.

### 4.1. Points stationnaires du système dérivé

Supposons qu'un point stationnaire de $(\Sigma)$ existe et notons-le $(a; b; c)$. C'est *a priori* un élément de $[0;1]^3$, qui est solution du système d'équations algébriques :

$$(\Sigma') : \begin{cases} a = 1 - bc \\ b = 1 - ac \\ c = 1 - ab \end{cases}$$

Supposons que l'un des trois réels est égal à 1. On montre alors facilement que $(1; 0; 1)$, $(1; 1; 0)$ et $(0; 1; 1)$ sont les seules solutions. Supposons désormais qu'aucun des trois réels n'est égal à 1. Aucun n'est égal à 0. Il s'agit donc de résoudre $(\Sigma')$ dans $]0;1[^3$. On a alors :

$$(\Sigma') \Leftrightarrow \begin{cases} a = 1 - bc \\ a + bc = b + ac \\ a + bc = c + ab \end{cases} \Leftrightarrow \begin{cases} a = 1 - bc \\ (a - b)(1 - c) = 0 \\ (a - c)(1 - b) = 0 \end{cases} \Leftrightarrow \begin{cases} a = b = c \\ a = 1 - a^2 \end{cases}$$

Les nombres $a$, $b$ et $c$ sont donc l'unique solution $\alpha_3$ dans $]0;1[$ de l'équation $(E_3) : x^2 + x - 1 = 0$ :

$$\alpha_3 = \frac{\sqrt{5} - 1}{2} = \varphi - 1 \, , \;\; \text{où} \;\; \varphi = \frac{1 + \sqrt{5}}{2} \;\; \text{est le nombre d'or.}$$

La réciproque est claire. Est ainsi établie l'existence de quatre points stationnaires pour $(\Sigma)$, seul le dernier étant dans $]0;1[^3$ : $(1; 1; 0), (0; 1; 1), (1; 0; 1)$ et $(\alpha; \alpha; \alpha)$, où on note ici $\alpha = \alpha_3$. D'où quatre points stationnaires pour $(S)$ : $(1; 0; 0), (0; 1; 0), (0; 0; 1)$ et $(1 - \alpha; 1 - \alpha; 1 - \alpha)$.

### 4.2. Dynamique du système dérivé au voisinage de ses points stationnaires

Étudions maintenant la stabilité du système $(\Sigma)$ au voisinage de ses points stationnaires.

Considérons le point $(\alpha; \alpha; \alpha)$, le seul dont la considération va s'avérer importante. Pour tout $m \in \mathbb{N}$, on notera désormais : $a_m = u_m - \alpha \;\; ; \;\; b_m = v_m - \alpha \;\; ; \;\; c_m = w_m - \alpha$. On a alors :

$$\forall m \in \mathbb{N}, \begin{cases} a_{m+1} + \alpha = 1 - (b_m + \alpha)(c_m + \alpha) \\ b_{m+1} + \alpha = 1 - (a_m + \alpha)(c_m + \alpha) \\ c_{m+1} + \alpha = 1 - (a_m + \alpha)(b_m + \alpha) \end{cases}$$

Le système récurrent linéarisé au voisinage de $(0; 0; 0)$ s'écrit donc, compte tenu de $1 - \alpha - \alpha^2 = 0$ :

$$\forall\, m \in \mathbb{N}, \quad \begin{cases} a_{m+1} = -\alpha b_m - \alpha c_m \\ b_{m+1} = -\alpha b_m - \alpha c_m \\ c_{m+1} = -\alpha b_m - \alpha c_m \end{cases}$$

La matrice de ce système est réelle et symétrique, donc diagonalisable dans $\mathbb{R}$. On obtient sans difficulté que son polynôme caractéristique est $\chi_M = (X - \alpha)^2 (X + 2\alpha)$. Elle a donc pour diagonalisée

$$D = \begin{pmatrix} -2\alpha & 0 & 0 \\ 0 & \alpha & 0 \\ 0 & 0 & \alpha \end{pmatrix}$$

Or, la suite de matrices de terme général $D^m = \begin{pmatrix} (-2\alpha)^m & 0 & 0 \\ 0 & \alpha^m & 0 \\ 0 & 0 & \alpha^m \end{pmatrix}$ diverge compte tenu du fait que $\alpha = \varphi - 1 \in\, ]0; 1[$ et que $-2\alpha = 2(1 - \varphi) < -1$. Donc le point stationnaire $(\alpha; \alpha; \alpha)$ est exponentiellement instable. Ces conclusions s'appliquent au point stationnaire de $(S)$ correspondant.

### 4.3. Convergence et limite de la suite duale

On peut toujours supposer, à symétrie près des rôles, que $0 < u_0 \leq v_0 \leq w_0 < 1$. On a alors :

$$\forall\, m \in \mathbb{N},\ 0 < u_m \leq v_m \leq w_m < 1$$

Cette propriété est initialisée. Supposons qu'il existe $m \in \mathbb{N}$ tel que $0 < u_m \leq v_m \leq w_m < 1$. On a :

$$u_{m+1} - v_{m+1} = (1 - v_m w_m) - (1 - u_m w_m) = w_m (u_m - v_m) \leq 0$$

Donc $u_{m+1} \leq v_{m+1}$. De plus, il est clair que $u_{m+1} = 1 - v_m w_m > 0$. On montre de même que $v_{m+1} \leq w_{m+1} < 1$. Ceci établit l'hérédité de la propriété considérée, qui est donc vraie par récurrence.

Observons de plus que, toujours avec les conventions précédentes, on a pour tout $m \in \mathbb{N}$ :

$$u_m v_{m+2} - v_m u_{m+2} = u_m v_m (u_m + w_m - u_m v_m w_m) - v_m u_m (v_m + w_m - u_m v_m w_m)$$
$$= u_m v_m (u_m - v_m) \leq 0$$

On en déduit, en tenant aussi compte de ce qui précède, que :

$$\forall\, m \in \mathbb{N},\ 1 \leq \frac{v_{m+2}}{u_{m+2}} \leq \frac{v_m}{u_m}$$

La suite $(v_{2q}/u_{2q})_{q \in \mathbb{N}}$, décroissante et minorée par 1, converge vers un réel $L \geq 1$. De même, les suites $(w_{2q}/v_{2q})_{q \in \mathbb{N}}$ et $(w_{2q}/u_{2q})_{q \in \mathbb{N}}$ convergent en décroissant vers des réels $L' \geq 1$ et $L'' \geq 1$.

Nous allons établir que $L = L' = L'' = 1$ en montrant au préalable le résultat suivant :

$$\forall\, m \in \mathbb{N},\ \frac{\frac{w_{m+2}}{u_{m+2}} - 1}{\frac{w_m}{u_m} - 1} < \frac{1}{2}$$

et pour tout $q \in \mathbb{N},\ 0 < \frac{w_{2q}}{u_{2q}} - 1 < \left(\frac{1}{2}\right)^q \left(\frac{w_0}{u_0} - 1\right)$

Posons, pour tout $m \in \mathbb{N}$, $K_m = v_m - u_m v_m w_m$ et $C_m = u_m w_m$. On a $K_m > 0$ et $C_m > 0$. On vérifie que : $w_{m+2} = K_m w_m + C_m$ et $u_{m+2} = K_m u_m + C_m$. On obtient alors, puisque l'irrégularité de la suite barypolygonale initiale implique l'existence du premier rapport (l'hypothèse d'une régularité à un rang donné impliquant la régularité initiale par récurrence descendante) :

$$\frac{\frac{w_{m+2}}{u_{m+2}} - 1}{\frac{w_m}{u_m} - 1} = \frac{K_m u_m}{K_m u_m + C_m} = 1 - \frac{C_m}{K_m u_m + C_m}$$

avec :

$$\frac{C_m}{K_m u_m + C_m} = \frac{u_m w_m}{u_m w_m + u_m v_m - u_m^2 v_m w_m} = \frac{w_m}{w_m + v_m - u_m v_m w_m} = \frac{1}{1 + \frac{v_m}{w_m} - u_m v_m}$$

Égalité qui implique, puisque $u_m v_m > 0$, que :

$$\frac{C_m}{K_m u_m + C_m} > \frac{1}{1 + \frac{v_m}{w_m}}$$

On a donc :

$$\left(\forall\, m \in \mathbb{N},\ \frac{v_m}{w_m} \leq 1\right) \Rightarrow \left(\forall\, m \in \mathbb{N},\ \frac{C_m}{K_m u_m + C_m} > \frac{1}{2}\right)$$

On en déduit bien :

$$\forall\, m \in \mathbb{N},\ 0 < \frac{w_{m+2}}{u_{m+2}} - 1 < \frac{1}{2}\left(\frac{w_m}{u_m} - 1\right)$$

Puis la seconde inégalité annoncée par récurrence immédiate.

D'où la convergence exponentielle vers 1 de la suite $(w_{2q}/u_{2q})_{q\in\mathbb{N}}$. Ainsi, $L'' = 1$. De plus :

$$\forall\, q \in \mathbb{N},\ 1 \leq \frac{w_{2q}}{v_{2q}} \leq \frac{w_{2q}}{u_{2q}}$$

On en déduit par passage à la limite que $L' = 1$. Il en résulte aussi immédiatement que $L = 1$. Avec dans les deux cas des convergences qui s'opèrent aussi bien de manière exponentielle.

On démontre de même que les suites $(w_{2q+1}/u_{2q+1})_{q\in\mathbb{N}}$, $(w_{2q+1}/v_{2q+1})_{q\in\mathbb{N}}$ et $(v_{2q+1}/u_{2q+1})_{q\in\mathbb{N}}$ convergent exponentiellement vers 1. On peut dès lors en déduire que :

$$\lim_{m\to+\infty} \frac{w_m}{u_m} = \lim_{m\to+\infty} \frac{w_m}{v_m} = \lim_{m\to+\infty} \frac{v_m}{u_m} = 1$$

Observons enfin que pour tout $q \in \mathbb{N}$ :

$$t_1^{(2q+1)} = \left(1 - t_2^{(2q)}\right)\left(1 - t_3^{(2q)}\right) = v_{2q}\, w_{2q}\ \text{et}\ t_2^{(2q+1)} = \left(1 - t_1^{(2q)}\right)\left(1 - t_3^{(2q)}\right) = u_{2q}\, w_{2q}$$

Ce qui implique :

$$\forall\, q \in \mathbb{N},\quad \frac{t_1^{(2q+1)}}{t_2^{(2q+1)}} = \frac{v_{2q}}{u_{2q}}$$

Et donc d'après ce qui précède :

$$\lim_{q\to+\infty} \frac{t_1^{(2q+1)}}{t_2^{(2q+1)}} = 1$$

Le même résultat peut être établi pour les rapports analogues. D'où résulte finalement :

$$\lim_{m\to+\infty} \frac{t_1^{(m)}}{t_2^{(m)}} = \lim_{m\to+\infty} \frac{t_2^{(m)}}{t_3^{(m)}} = \lim_{m\to+\infty} \frac{t_3^{(m)}}{t_1^{(m)}} = 1$$

On a ainsi montré que la duale de la suite barypolygonale $\mathfrak{B}$ converge exponentiellement vers le centre de gravité $G$ de $\mathcal{A}$. Ce résultat va permettre de préciser le comportement du système $(\Sigma)$.

## 4.4. Lieu des éléments du système dérivé

Si les $p$ réels du terme initial $t$ sont égaux à $1 - \alpha$, la suite $t$-barypolygonale initialisant la suite de ses dérivées est régulière : on a vu qu'alors $(S)$ est stationnaire en $(1 - \alpha; 1 - \alpha; 1 - \alpha)$. Étudions dans le cas contraire le comportement de $(\Sigma)$ pour une position donnée des termes d'un triplet $(u_m; v_m; w_m)$ par rapport au point « répulsif » $(\alpha; \alpha; \alpha)$, en commençant par déterminer leur lieu.

Supposons d'abord qu'il existe $m_0 \in \mathbb{N}$ tel que $(u_{m_0}; v_{m_0}; w_{m_0}) \in {]\alpha; 1[}^3$. On montre alors par récurrence (en présentant ici deux formes de justification) que pour tout $q \in \mathbb{N}$ :

$$(u_{m_0+2q}; v_{m_0+2q}; w_{m_0+2q}) \in {]\alpha; 1[}^3 \text{ et } (u_{m_0+2q+1}; v_{m_0+2q+1}; w_{m_0+2q+1}) \in {]0; \alpha[}^3$$

Cette propriété est initialisée. Compte tenu de $1 - \alpha^2 = \alpha$, on a en effet :

$$u_{m_0+1} = 1 - v_{m_0} w_{m_0} = 1 - \big((v_{m_0} - \alpha) + \alpha\big)\big((w_{m_0} - \alpha) + \alpha\big)$$
$$= \alpha - \alpha(v_{m_0} - \alpha) - \alpha(w_{m_0} - \alpha) - (v_{m_0} - \alpha)(w_{m_0} - \alpha)$$

Il est alors clair, compte tenu de $(u_{m_0}; v_{m_0}; w_{m_0}) \in {]\alpha; 1[}^3$ et de $\alpha > 0$, que $u_{m_0+1} < \alpha$. De même, $v_{m_0+1} < \alpha$ et $w_{m_0+1} < \alpha$, de sorte que $(u_{m_0+1}; v_{m_0+1}; w_{m_0+1}) \in {]0; \alpha[}^3$. Soit maintenant $q \in \mathbb{N}$ tel que $(u_{m_0+2q}; v_{m_0+2q}; w_{m_0+2q}) \in {]\alpha; 1[}^3$. On établit comme précédemment que $(u_{m_0+2q+1}; v_{m_0+2q+1}; w_{m_0+2q+1}) \in {]0; \alpha[}^3$. De plus :

$$u_{m_0+2q+2} = 1 - v_{m_0+2q+1} w_{m_0+2q+1}$$

Or : $0 < v_{m_0+2q+1} < \alpha$ et $0 < w_{m_0+2q+1} < \alpha$ impliquent $0 < v_{m_0+2q+1} w_{m_0+2q+1} < \alpha^2$ et donc :

$$\alpha = 1 - \alpha^2 < 1 - v_{m_0+2q+1} w_{m_0+2q+1} < 1$$

Par conséquent : $u_{m_0+2q+2} \in {]\alpha; 1[}$. De même, $v_{m_0+2q+2} \in {]\alpha; 1[}$ et $w_{m_0+2q+2} \in {]\alpha; 1[}$. On peut alors de nouveau établir comme dans le cas où $q = 0$ que $(u_{m_0+2q+3}; v_{m_0+2q+3}; w_{m_0+2q+3}) \in {]0; \alpha[}^3$.

On établit bien sûr de manière similaire par récurrence que s'il existe $m_0 \in \mathbb{N}$ tel que $(u_{m_0}; v_{m_0}; w_{m_0}) \in {]0; \alpha[}^3$, alors pour tout $q \in \mathbb{N}$ :

$$(u_{m_0+2q}; v_{m_0+2q}; w_{m_0+2q}) \in {]0; \alpha[}^3 \text{ et } (u_{m_0+2q+1}; v_{m_0+2q+1}; w_{m_0+2q+1}) \in {]\alpha; 1[}^3$$

Reste alors à savoir si l'existence d'un entier $m_0$ satisfaisant l'une des deux hypothèses examinées est bien toujours assurée. Montrons que la réponse est positive.

On reprend l'hypothèse $0 < u_0 \leq v_0 \leq w_0 < 1$. D'après le théorème de Bolzano-Weierstrass, la suite $(u_m; v_m; w_m)_{m \in \mathbb{N}}$ d'éléments de $[0; 1]^3$ admet une suite extraite $\big((u_{\psi(m)}); (v_{\psi(m)}); (w_{\psi(m)})\big)$ convergeant vers un triplet $(a; b; c)$ de $[0; 1]^3$, ici tel que $0 \leq a \leq b \leq c \leq 1$. Comme on a supposé $t$ distinct de $(1 - \alpha; 1 - \alpha; 1 - \alpha)$ et puisque le point stationnaire $(\alpha; \alpha; \alpha)$ est exponentiellement instable, on a $(a; b; c) \neq (\alpha; \alpha; \alpha)$. Supposons que $a < \alpha$. Il résulte du 4.3 que :

$$\lim_{m \to +\infty} \frac{u_{\psi(m)}}{v_{\psi(m)}} = 1 = \lim_{m \to +\infty} \frac{v_{\psi(m)}}{w_{\psi(m)}}$$

On en déduit que si $a = 0$, alors $b = c = 0$. Et que si $a > 0$, alors on a encore $a = b = c$. Dans les deux cas, on peut donc écrire :

$$\forall \, \varepsilon > 0, \ \exists \, M \in \mathbb{N}, \ \forall \, m \geq M, \ 0 < u_{\psi(m)} \leq v_{\psi(m)} \leq w_{\psi(m)} < a + \varepsilon$$

En choisissant $\varepsilon = \frac{\alpha - a}{2}$, on obtient ainsi :

$$\forall\, m \geq M,\ 0 < u_{\psi(m)} \leq v_{\psi(m)} \leq w_{\psi(m)} < \frac{a + \alpha}{2} < \alpha$$

L'entier $m_0 = \psi(M)$ satisfait alors la condition recherchée. Un tel entier existe aussi lorsque $a \geq \alpha$, $m_0$ étant cette fois tel que $(u_{m_0}; v_{m_0}; w_{m_0}) \in\, ]\alpha; 1[^3$ : il suffit d'utiliser le fait que $\alpha < c \leq 1$ et d'en déduire des inégalités analogues aux précédentes sur l'intervalle $]\alpha; 1[$.

On peut alors conclure de l'ensemble des considérations précédentes qu'il existe bien dans tous les cas $m_0 \in \mathbb{N}$ tel que pour tout $q \in \mathbb{N}$ :

$$\left((u_{m_0+2q}; v_{m_0+2q}; w_{m_0+2q}) \in\, ]0; \alpha[^3\ \text{et}\ (u_{m_0+2q+1}; v_{m_0+2q+1}; w_{m_0+2q+1}) \in\, ]\alpha; 1[^3\right)$$
$$\text{ou}$$
$$\left((u_{m_0+2q}; v_{m_0+2q}; w_{m_0+2q}) \in\, ]\alpha; 1[^3\ \text{et}\ (u_{m_0+2q+1}; v_{m_0+2q+1}; w_{m_0+2q+1}) \in\, ]0; \alpha[^3\right)$$

Ces conclusions se transfèrent bien sûr par « conjugaison » au système dérivé $(S)$.

## 4.5. Convergence des suites extraites des termes de rang pair et de rang impair d'éléments du système dérivé

On peut désormais étudier la convergence des suites extraites des termes de rang pair et de rang impair d'éléments du système dérivé.

Les discussions précédentes montrent qu'il suffit, sans perte de généralité, d'examiner le cas d'un entier $m_0 \in \mathbb{N}$ tel que $0 < u_{m_0} \leq v_{m_0} \leq w_{m_0} < \alpha$, ce qui implique que $u_m \leq v_m \leq w_m$ pour tout $m \geq m_0$. Montrons alors que les trois suites $(u_{m_0+2q})_{q \in \mathbb{N}}$, $(v_{m_0+2q})_{q \in \mathbb{N}}$ et $(w_{m_0+2q})_{q \in \mathbb{N}}$ convergent vers 0, et que les trois suites $(u_{m_0+2q+1})_{q \in \mathbb{N}}$, $(v_{m_0+2q+1})_{q \in \mathbb{N}}$ et $(w_{m_0+2q+1})_{q \in \mathbb{N}}$ convergent vers 1.

Considérons $f : x \mapsto 1 - x^2$, bijection décroissante de $[0; 1]$ dans $[0; 1]$. Si $u_{m_0+1} > f(w_{m_0})$, on pose $\tau = w_{m_0}$ ; sinon, on pose $\tau = f^{-1}(u_{m_0+1})$. On définit alors la suite $(\tau_m)_{m \geq m_0}$ par :

$$\begin{cases} \tau_{m_0} = \tau \\ \forall\, m \geq m_0,\ \tau_{m+1} = f(\tau_m) \end{cases}$$

Montrons par récurrence que :

Pour tout $q \in \mathbb{N}$, $\tau_{m_0+2q} \geq w_{m_0+2q}$ et $\tau_{m_0+2q+1} \leq u_{m_0+2q+1}$.

Ces deux inégalités sont vraies pour $q = 0$, par définition de $\tau$. La première l'est évidemment. Quant à la seconde, elle découle du fait que

$$\left(u_{m_0+1} > f(w_{m_0})\right) \Rightarrow \left(u_{m_0+1} > f(\tau) = f(\tau_{m_0}) = \tau_{m_0+1}\right)$$

$$\text{et}\ \left(u_{m_0+1} \leq f(w_{m_0})\right) \Rightarrow (u_{m_0+1} = f(\tau) = \tau_{m_0+1})$$

De plus, si on considère un entier $q \in \mathbb{N}$ tel que $\tau_{m_0+2q} \geq w_{m_0+2q}$, alors :

$$u_{m_0+2q+1} = f\left(\sqrt{v_{m_0+2q}\, w_{m_0+2q}}\right) \geq f\left(w_{m_0+2q}\right) \geq f(\tau_{m_0+2q}) = \tau_{m_0+2q+1}$$

et obtient de même que $\tau_{m_0+2q+2} \geq w_{m_0+2q+2}$, puis que $\tau_{m_0+2q+3} \leq u_{m_0+2q+3}$. Donc la propriété annoncée est aussi héréditaire.

Or, d'après le théorème 3, on sait que la suite $(\tau_{2m})$ converge vers 0 et que la suite $(\tau_{2m+1})$ converge vers 1. Donc par comparaison, on déduit de ce qui précède que :

$$\lim_{m \to +\infty} u_{2m+1} = 1 \text{ et } \lim_{m \to +\infty} w_{2m} = 0$$

Ce qui implique, là encore par comparaison, qu'on a aussi :

$$\lim_{m \to +\infty} v_{2m+1} = 1 = \lim_{m \to +\infty} w_{2m+1} \text{ et } \lim_{m \to +\infty} u_{2m} = 0 = \lim_{m \to +\infty} v_{2m}$$

Les résultats analogues (mais intervertis dans les rôles) se montrent bien sûr si l'on suppose que $m_0$ est tel que $\alpha < u_{m_0} \leq v_{m_0} \leq w_{m_0} < 1$. Le comportement de $(\Sigma)$ est ainsi déterminé dans tous les cas.

On déduit par « conjugaison » des résultats obtenus aux 4.4 et 4.5 que les termes du système dérivé $(S)$ satisfont les mêmes propriétés que celles établies pour $(\Sigma)$, à inversion près des rôles.

## 5. Conclusion

Les problèmes posés au 3. ont ainsi été entièrement résolus lorsque $p = 3$, les conjectures alors émises ayant été corroborées. La seconde partie de cet article montrera que ces résultats sont conservés lorsque $p \geq 4$, en généralisant la méthode employée dans la précédente étude.

## 6. Bibliographie